\documentstyle[amscd,amssymb,verbatim]{amsart}
\newcommand{\lrar}[1]{\begin{picture}(50,10)(-25,-5)                          
\put(-25,0){\vector(1,0){50}}
\put(0,5){\makebox(0,0)[b]{\mbox{$#1$}}}
\end{picture}}

\newcommand{\ldar}[1]{\begin{picture}(10,50)(-5,-25)
\put(0,25){\vector(0,-1){50}}
\put(5,0){\mbox{$#1$}}
\end{picture}}

\newcommand{\luar}[1]{\begin{picture}(10,50)(-5,-25)
\put(0,-25){\vector(0,1){50}}
\put(5,0){\mbox{$#1$}}
\end{picture}}

\newcommand{\Ob}{\operatorname{Ob}}

\newcommand{\ti}{\tilde}

\newcommand{\Vect}{\operatorname{Vect}}

\newcommand{\CC}{{\cal C}}

\newcommand{\Tr}{\operatorname{Tr}}

\newcommand{\Pic}{\operatorname{Pic}}

\newcommand{\Id}{\operatorname{Id}}
\newcommand{\ga}{\gamma}

\newcommand{\eps}{\epsilon}

\numberwithin{equation}{section}

\newtheorem{thm}{Theorem}[section]
\newtheorem{prop}[thm]{Proposition}
\newtheorem{lem}[thm]{Lemma}

\newenvironment{rem}{\vspace{3mm}\noindent
{\bf Remark.}}{\vspace{3mm}}

\theoremstyle{definition}
\newtheorem{defi}[thm]{Definition}
\newcommand{\Pf}{\noindent {\it Proof}}
\newcommand{\id}{\operatorname{id}}

\newcommand{\ov}{\overline}
\newcommand{\we}{\wedge}
\renewcommand{\Im}{\operatorname{Im}}

\newcommand{\ra}{\rightarrow}

\newcommand{\VV}{{\cal V}}

\newcommand{\LL}{{\cal L}}
\newcommand{\UU}{{\cal U}}
\renewcommand{\O}{{\cal O}}
\newcommand{\Om}{\Omega}

\newcommand{\pa}{\partial}
\newcommand{\Hom}{\operatorname{Hom}}
\newcommand{\Ext}{\operatorname{Ext}}

\renewcommand{\a}{\alpha}
\renewcommand{\b}{\beta}
\newcommand{\om}{\omega}

\newcommand{\la}{\lambda}
\newcommand{\th}{\theta}
\newcommand{\C}{{\Bbb C}}
\newcommand{\R}{{\Bbb R}}
\newcommand{\Z}{{\Bbb Z}}

\newcommand{\wt}{\widetilde}

\newcommand{\sign}{\operatorname{sign}}

\newcommand{\ed}{\qed\vspace{3mm}}
\newcommand{\Ax}{\operatorname{Ax}}

\title{$A_{\infty}$-structures on an elliptic curve}
\author{A. Polishchuk}
\thanks{This work is partially supported by NSF grant}

\begin{document}

\begin{abstract}
The main result of this paper is the proof of the "transversal part"
of the homological mirror symmetry conjecture for an elliptic curve
which states an equivalence of two $A_{\infty}$-structures on the
category of vector bundles on an elliptic curves. The proof is
based on the study of $A_{\infty}$-structures on the category of line
bundles over an elliptic curve satisfying some natural restrictions 
(in particular, $m_1$ should be zero, $m_2$ should coincide with the 
usual composition). The key observation is that such a structure is uniquely 
determined up to homotopy by certain triple products. 
\end{abstract}

\maketitle

\bigskip

\centerline{\sc Introduction}

\bigskip

Let $E$ be an elliptic curve over a field $k$. Let us denote by $\Vect(E)$
the category of algebraic vector bundles on $E$, where as space of morphisms
from $V_1$ to $V_2$ we take the graded space $\Hom(V_1,V_2)\oplus\Ext^1(V_1,V_2)$
with the natural composition law. In this paper we study extensions of
this (strictly associative) composition to $A_{\infty}$-structures on $\Vect(E)$
(see section \ref{ainf} for the definition).
The motivation comes from the homological mirror symmetry for 
elliptic curves formulated by Kontsevich (see \cite{K})
which provides two such extensions in the case $k=\C$
and states that they should be equivalent.
We recall the definitions of these $A_{\infty}$-structures in section 
\ref{two}.
One of them is an $A_{\infty}$-version of the derived category
of vector bundles while another comes from a general construction
in symplectic geometry due to Fukaya. 
Roughly speaking, one can associate to an indecomposable 
vector bundle on
$E$ a geodesic circle on the torus $\R^2/\Z^2$ with a local
system on it. Then the second $A_{\infty}$-structure is defined using
generating series counting holomorphic maps from the disk bounding given 
geodesic circles.

In \cite{PZ} we checked that the double products defined in this way coincide with
the standard composition law on $\Vect(E)$.
In this paper we use this together with some calculations of triple products
(see \cite{P-Mas}) to
prove the essential part of the homological mirror conjecture for $E$. Namely,
we construct a homotopy between {\it transversal} products given by these
two $A_{\infty}$-structures. This means that we are looking only at the products
such that the corresponding geodesic circles form a transversal configuration.
The advantage is that in this case the homotopy can be constructed in a canonical 
way. We leave to a future investigation more subtle points of defining
non-transversal products in the Fukaya category and extending the above homotopy
to the entire derived categories.

Note that the equality of double products and triple Massey products in 
$A_{\infty}$-categories corresponding to a mirror dual pair
(symplectic torus, abelian variety) was established by Fukaya in \cite{F3}.
In the case of elliptic curves, as we show in the present paper,  
this is enough for (transversal part of) the
homological mirror conjecture. For abelian varieties of higher
dimensions, a version of this conjecture 
was recently proved by Kontsevich and Soibelman in \cite{KS}.
\footnote{The $A_{\infty}$-equivalence established in \cite{KS}
deals with certain full subcategories in symplectic and holomorphic
$A_{\infty}$-categories. In fact, both sides are slightly
modified: Fukaya category is replaced by its degeneration, while
on the holomorphic side the ground field is changed from $\C$ to
$\C((q))$. Also, only transversal products are considered.}
The main point of our paper is that in the case of elliptic curves
we can formulate a result on $A_{\infty}$-structures on the category $\LL$ 
of line bundles on $E$ which is valid over an arbitrary field $k$.
More precisely, 
we axiomatize the notion of transversality and prove that if one imposes some
natural restrictions on a transversal $A_{\infty}$-structure on $\LL$ (in particular
$m_1=0$, $m_2$ is equal to the standard composition), 
then such a structure is uniquely
determined (up to homotopy) by certain triple products. Namely, these are triple
products which are invariant under any homotopy.  
We apply this result to two $A_{\infty}$-structures
arising in the homological mirror symmetry and then use the isogenies between elliptic
curves (as in \cite{PZ}) to construct the required homotopy on the category of 
vector bundles on $E$.

The natural framework for the generalization of our result which is valid over any 
field $k$ should
involve the notion of a triangulated $A_{\infty}$-category (as sketched in \cite{K-t}). 
Our result seems to imply that there exists a unique up to homotopy triangulated 
$A_{\infty}$-structure on the derived category of an elliptic curve which is 
compatible with the standard products and with Serre duality 
(see section \ref{cyc-sec} for the definition of the latter compatibility). 
Indeed, the triple products appearing in our statement
are univalued Massey products which are uniquely determined by the double products
in the case when $A_{\infty}$-structure is triangulated. 
One may hope that such a uniqueness of $A_{\infty}$-structure on the derived category
holds for other varieties (e.g. for abelian varieties of arbitrary dimension).
The main reason why in the case of $A_{\infty}$-structures on elliptic 
curve only triple products matter is the absense of non-trivial
univalued well-defined $k$-tuple Massey products for $k>3$
\footnote{In fact, there exist non-zero univalued quadruple Massey products on 
elliptic curve but they can be expressed via triple products.}. 

Conventions: We always work over a ground field $k$; we 
specialize to $k=\C$ when talking about homological mirror symmetry.
To shorten the formulas sometimes we denote 
the tensor product of vector spaces $V_1$ and $V_2$ over $k$
simply by $V_1V_2$ omitting the sign of the tensor product. We use
the same abbreviation for tensor products of vector bundles.
By a {\it bundle} we always mean an algebraic vector bundle 
(or a holomorphic vector bundle
if $k=\C$). When working with $A_{\infty}$-categories it is convenient
to denote the $n$-tuple products of composable morphisms
$a_1:X_0\ra X_1$, $a_2:X_1\ra X_2$, ..., $a_n:X_{n-1}\ra X_n$ by
$m_n(a_1,a_2,\ldots,a_n)$. In particular, we denote the double composition
by $m_2(a_1,a_2)$ which we often abbreviate to $a_1a_2$. This contradicts
to the usual convention of going from right to left when considering
composition in the usual categories. To avoid confusion
we will use the notation $a_2\circ a_1$ for the composition in
the usual categories. 

\section{$A_{\infty}$-structures and their homotopies}
\label{ainf}

In the following definitions we use the sign convention of \cite{GJ} which is
different from the one in the original definition of \cite{S}.

\subsection{$A_{\infty}$-algebras}
A ($\Z$-graded) $A_{\infty}$-algebra is a $\Z$-graded
vector space $A$ 
equipped with linear maps $m_k: A^{\otimes k}\ra A$ for
$k\ge 1$ of degree $2-k$ satisfying for every $n\ge 1$
the following $A_{\infty}$-constraint $\Ax_n$:
$$\sum_{k+l=n+1}\sum_{j=1}^{k}
(-1)^{l(\ti{a}_1+\ldots+\ti{a}_{j-1})+j(l+1)}
m_k(a_1,\ldots,a_{j-1},m_l(a_j,\ldots,a_{j+l-1}),a_{j+l},\ldots,a_n)=0$$
where $\ti{a}_i=\deg(a_i)\mod(2)$. For example, $\Ax_1$ says that $m_1^2=0$,
$\Ax_2$ gives the Leibnitz identity for $m_1$ and $m_2$, etc.
One can consider $m_n$
as components of a coderivation $d$ of the coalgebra $T(sA)$ where
$s$ denote the suspension. The elements of $T(sA)$ are denoted traditionally
as follows: 
$$[a_1|a_2|\ldots|a_k]=(sa_1)\otimes (sa_2)\otimes\ldots (sa_k).$$
The coderivation $d$ has a component $d_k:(sA)^{\otimes k}\ra sA$ given by
$s\circ m_k\circ (s^{-1})^{\otimes k}$, so that
$$d([a_1|\ldots|a_n])= \sum_{k+l=n+1}\sum_{j=1}^{k}
(-1)^{\ti{a}_1+\ldots+\ti{a}_{j-1}+j-1+\mu(a_j,\ldots,a_{j+l-1})}
[a_1|\ldots|a_{j-1}|m_l(a_j,\ldots,a_{j+l-1})|a_{j+l}|\ldots|a_n].$$
where for every collection of elements $(a_1,\ldots,a_k)$ in $A$
we denote
$$\mu(a_1,\ldots,a_k)=(k-1)\ti{a}_1+(k-2)\ti{a}_2+\ldots +\ti{a}_{k-1}+\frac{k(k-1)}{2}.$$
The $A_{\infty}$-constraints are equivalent to the condition $d^2=0$.

For a pair of $A_{\infty}$-algebras $(A,m^A)$ and $(B,m^B)$ there is a natural
notion of a $A_{\infty}$-morphism from $A$ to $B$. Namely,
such a morphism consists of the data $(f_n, n\ge 1)$ where
$f_n:A^{\otimes n}\ra B$ is a linear map of degree $1-n$ such that
\begin{align*}
&\sum_{1\le k_1<k_2<\ldots<k_i=n}(-1)^{\eps_L}
m_i^B(f_{k_1}(a_1,\ldots,a_{k_1}),f_{k_2-k_1}(a_{k_1+1},\ldots,a_{k_2}),
\ldots,f_{n-k_{i-1}}(a_{k_{i-1}+1},\ldots,a_n))\\
&=\sum_{k+l=n+1}\sum_{j=1}^{k}(-1)^{\eps_R}
f_k(a_1,\ldots,a_{j-1},m_l^A(a_j,\ldots,a_{j+l-1}),a_{j+l},\ldots,a_n)
\end{align*}
where the signs $\eps_L$ and $\eps_R$ are defined as follows:
\begin{align*}
&\eps_L=\mu(a_1,\ldots,a_{k_1})+\mu(a_{k_1+1},\ldots,a_{k_2})+\ldots
+\mu(a_{k_{i-1}+1},\ldots,a_n)+\\
&\mu(f_{k_1}(a_1,\ldots,a_{k_1}),
\ldots,f_{n-k_{i-1}}(a_{k_{i-1}+1},\ldots,a_n)),
\end{align*}
$$\eps_R=\wt{a}_1+\ldots+\wt{a}_{j-1}+j-1+\mu(a_j,\ldots,a_{j+l-1})+
\mu(a_1,\ldots,a_{j-1},m_l^A(a_j,\ldots,a_{j+l-1}),a_{j+l},\ldots,a_n).$$
Again one can consider $(f_n)$ as components of a coalgebra homomorphism 
$F:T(sA)\ra T(sB)$, so that the above equation is equivalent to 
$$F\circ d^A=d^B\circ F,$$ 
where $d_A$ (resp. $d^B$) is the coderivation on $A$ (resp. $B$)
defined by $m^A$ (resp. $m^B$).
In particular, there is a natural
composition of $A_{\infty}$-morphisms defined as follows:
\begin{align*}
&(f\circ g)_n(a_1,\ldots,a_n)=\\
&\sum_{1\le k_1<k_2<\ldots<k_i=n}(-1)^{\eps}
f_i(g_{k_1}(a_1,\ldots,a_{k_1}),g_{k_2-k_1}(a_{k_1+1},\ldots,a_{k_2}),
\ldots,g_{n-k_{i-1}}(a_{k_{i-1}+1},\ldots,a_n))
\end{align*}
where
\begin{align*}
&\eps=\mu(a_1,\ldots,a_{k_1})+\mu(a_{k_1+1},\ldots,a_{k_2})+\ldots+
\mu(a_{k_{i-1}+1},\ldots,a_n)+\\
&\mu(g_{k_1}(a_1,\ldots,a_{k_1}),\ldots,
g_{n-k_{i-1}}(a_{k_{i-1}+1},\ldots,a_n)).
\end{align*}

In the case when $B$ and $A$ have the same
underlying spaces and $f_1=\id$ we will call
the data $f=(f_n,n\ge 2)$ a {\it homotopy} between two
$A_{\infty}$-structures $m=m^A$ and $m'=m^B$ on the same space.
Note that for homotopic $m$ and $m'$ we necessarily have $m_1=m'_1$.
If $f$ is a homotopy between $m$ and $m'$, $g$ is a homotopy
between $m'$ and $m''$ then $g\circ f$ is a homotopy between $m$ and $m''$.

\begin{lem}\label{homotopy} Let
$m=(m_n)$ be an $A_{\infty}$-structure on $A$,
$(f_n:A^{\otimes n}\ra A,n\ge 2)$
be an arbitrary family of maps, $\deg f_n=1-n$. Then there exists
a unique $A_{\infty}$-structure $m'$ on $A$ such that $f=(f_n)$ (where
$f_1=\id$)
is a homotopy between $m$ and $m'$.
\end{lem}

\Pf . This follows immediately from the fact that the coalgebra homomorphism
$T(sA)\ra T(sA)$ defined by $(f_n)$
is an isomophism. 
\ed

We denote the $A_{\infty}$-structure $m'$ constructed
in the above lemma by $m+\delta f$ (note that it depends non-linearly on $f$).

\subsection{$A_{\infty}$-categories}
The definition of an $A_{\infty}$-category is similar to that
of an $A_{\infty}$-algebra (see \cite{F1}, \cite{Ke}). 
Namely, an $A_{\infty}$-category $\CC$
consists of a class of objects $\Ob\CC$,
for every pair of objects $X$ and $X'$ a graded space of morphisms
$\Hom(X,X')$, and a collection of linear maps (compositions)
$$m_k:\Hom(X_0,X_1)\otimes\Hom(X_1,X_2)\otimes\ldots\otimes
\Hom^*(X_{k-1},X_{k})\ra\Hom(X_0,X_{k})$$
of degree $2-k$ for all $k\ge 1$. The associativity constraint is
that these compositions define a structure of $A_{\infty}$-algebra
on $\oplus_{ij}\Hom(X_i,X_j)$ for every collection
$X_0,\ldots,X_n\in\Ob\CC$.

An $A_{\infty}$-functor (see \cite{F2}, \cite{Ke}) $\phi:\CC\ra\CC'$ between
$A_{\infty}$-categories consists of a map $\phi:\Ob\CC\ra\Ob\CC'$ and
of a collection of linear maps
$$f_k:\Hom_{\CC}(X_0,X_1)\otimes\Hom_{\CC}(X_1,X_2)\otimes\ldots\otimes
\Hom_{\CC}(X_{k-1},X_{k})\ra\Hom_{\CC'}(\phi(X_0),\phi(X_{k}))$$
of degree $1-k$ for $k\ge 1$,
which define $A_{\infty}$-morphisms
$\oplus_{ij}\Hom_{\CC}(X_i,X_j)\ra\oplus_{ij}\Hom_{\CC'}(\phi(X_i),
\phi(X_j))$.

Now assume that we are given two structures of $A_{\infty}$-category
with the same class of objects $\CC$ and with the same morphism spaces.
Let $m=(m_n)$ and $m'=m'_n$ be the collections of the corresponding composition maps. 
A homotopy between $m$ and $m'$ is an $A_{\infty}$-functor
$\phi:(\CC,m)\ra (\CC,m')$ such that the corresponding map on objects is identity
and such that $f_1$ is the identity map on morphisms.
The analogue of lemma \ref{homotopy} is valid in this situation.

In the case when $m_1=0$ for an $A_{\infty}$-category $\CC$ the products $m_2$
define a structure of the usual category on $\CC$ (without units).
If we have two such $A_{\infty}$-categories $\CC$ and $\CC'$ and a functor
$\phi_0:\CC\ra\CC'$ between them considered as usual categories then we say
that $\phi_0$ is {\it strictly compatible} with $A_{\infty}$-structures if
it extends to an $A_{\infty}$-functor with $f_k=0$ for $k>1$.

Let $X$ be an object of an $A_{\infty}$-category which has $m_1=0$
equipped with a decomposition
$X=X_1\oplus X_2$ into a direct sum. By definition (here we deal with the usual 
category structure without units) this means that
we have functorial in $Y$ isomorphisms
$$\Hom(X,Y)\simeq\Hom(X_1,Y)\oplus\Hom(X_2,Y)$$
and 
$$\Hom(Y,X)\simeq\Hom(Y,X_1)\oplus\Hom(Y,X_2).$$
We say that the decomposition $X=X_1\oplus X_2$ is
strictly compatibile with an $A_{\infty}$-structure
if every composition $m_n$ involving the spaces $\Hom(X,Y)$ or $\Hom(Y,X)$
is a direct sum of the corresponding compositions with the spaces $\Hom(X_i,Y)$
and $\Hom(Y,X_i)$.

\subsection{Cyclic $A_{\infty}$-structures}
\label{cyc-sec}

We will consider a special class of $A_{\infty}$-algebras,
namely, those equipped with a cyclic symmetry.

\begin{defi}
Let $A$ be an $A_{\infty}$-algebra equipped with a
bilinear form $b:A\otimes A\ra k$. 
We will call $A$ {\it cyclic} if for every $n\ge 1$
the following identity is satisfied:
\begin{equation}\label{cyclic}
b(m_n(a_1,\ldots,a_n),a_{n+1})=(-1)^{n(\ti{a}_1+1)}
b(a_1, m_n(a_2,\ldots,a_{n+1})).
\end{equation}
\end{defi}

\begin{rem}
Assume in addition that $b$ satisfies the following symmetry:
$$b(a_1,a_2)=(-1)^{\wt{a_1}\wt{a_2}}b(a_2,a_1).$$
Then (\ref{cyclic}) can be rewritten as follows:
$$b(m_n(a_1,\ldots,a_n),a_{n+1})=(-1)^{n+\wt{a_1}(\wt{a_2}+\ldots+\wt{a_{n+1}})}
b(m_n(a_2,\ldots,a_{n+1}),a_1).$$
\end{rem}

Using $b$ we can define a linear functional $\xi$ on $T(sA)$ by setting
$\xi=b\circ (s^{-1})^{\otimes 2}$ on $(sA)^{\otimes 2}$
while $\xi=0$ on $(sA)^{\otimes n}$ for $n\neq 2$. Thus, we have
$$\xi([a_1|a_2])=(-1)^{\wt{a}_1+1}b(a_1,a_2),$$
Then the equation (\ref{cyclic}) is equivalent to the condition
$\xi\circ d=0$ where $d$ is the coderivation defined by $(m_n)$.

The collection of maps $f=(f_n:A^{\otimes n}\ra A, n\ge 1)$,
$\deg f_n=1-n$, $f_1=\id$ is called a {\it cyclic homotopy} if 
\begin{equation}\label{cyclichom}
\sum_{k+l=n}(-1)^{(l+1)(\wt{a}_1+\ldots+\wt{a}_k)+nk} 
b(f_k(a_1,\ldots,a_k),f_{l}(a_{k+1},\ldots,a_n))=0
\end{equation}
for $n\ge 3$. This is equivalent to the condition $\xi\circ F=\xi$ where
$F:T(sA)\ra T(sA)$ is the coalgebra homomorphism defined by $(f_n)$.
Let $m$ be a cyclic $A_{\infty}$-structure, $f$ be a cyclic homotopy. 
Then the $A_{\infty}$-structure $m+\delta f$ is cyclic.
If $f$ and $g$ are cyclic homotopies then $f\circ g$ is also cyclic.

\begin{rem} Assume that $f_k=0$ unless $k=1$ or $k=n$ for some $n\ge 2$.
Then $f$ is a cyclic homotopy if and only if
$$b(f_n(a_1,\ldots,a_n),a_{n+1})=(-1)^{(n+1)\wt{a}_1+n}b(a_1,f_n(a_2,\ldots,a_{n+1}))$$
and 
$$b(f_n(a_1,\ldots,a_n),f_n(a_{n+1},\ldots,a_{2n}))=0.$$
\end{rem}

The definition of cyclic $A_{\infty}$-categories follows the same pattern.
We assume that there is a bilinear form
$$b:\Hom(X,Y)\otimes\Hom(Y,X)\ra k$$
for every pair of objects $(X,Y)$. Then an $A_{\infty}$-category is called cyclic
if the identity (\ref{cyclic}) is satisfied whenever 
$a_1\in\Hom(X_1,X_2),\ldots, a_n\in\Hom(X_n,X_{n+1})$, $a_{n+1}\in\Hom(X_{n+1},X_1)$.
Similarly we define cyclic homotopy between two cyclic $A_{\infty}$-structures
with the same objects and morphism spaces (and the same bilinear form $b$).

\subsection{Transversal $A_{\infty}$-structures}
Assume that we are given a class of objects
and a notion of transversality for pairs of objects.
We will call an $n$-tuple of objects $(X_1,\ldots,X_n)$ {\it transversal}
if for every $1\le i<j\le n$ the pair $(X_i,X_j)$ is transversal.
Then the structure of {\it transversal} (cyclic) $A_{\infty}$-category
on this class of objects consists of the following data. For every pair
of transversal objects $(X,Y)$ a graded space of morphisms
$\Hom(X,Y)$ is given. For every transversal collection $(X_0,\ldots,X_{n})$,
$n\ge 1$ we have linear maps
$$m_n:\Hom(X_0,X_1)\Hom(X_1,X_2)\ldots \Hom(X_{n-1},X_{n})\ra\Hom(X_0,X_{n})$$
of degree $2-n$ such that the axioms $\Ax_n$ 
and the identity (\ref{cyclic}) are satisfied whenever the objects involved
in it form a transversal collection.
Similarly we define a notion of homotopy between transversal 
$A_{\infty}$-structures and the cyclic analogues of these notions.

The motivating example is that of Fukaya category (see \cite{F1}) where
objects are Lagrangain submanifolds in a symplectic manifold with some
additional structure. Then we have the standard notion of transversality for pairs
of Lagrangians. However, notice that the notion of transversality for $n$-tuples
we use is weaker than the standard one: we just require every pair
of them to intersect transversally but, for example, we allow three Lagrangians
to intersect in one point. 

\subsection{Two $A_{\infty}$-structures on $\Vect(E)$}\label{two}

The first $A_{\infty}$-structure (or rather a class of equivalent
structures) can be defined on the category $\Vect(M)$ where $M$ is a variety
over $k$ or a complex manifold as follows.
Let us choose some functorial acyclic resolution $V\ra R^{\cdot}(V)$ for every 
vector bundle $V$ on $M$ such that for every pair of bundles
there are functorial morphisms
$$R^{\cdot}(V_1)\otimes R^{\cdot}(V_2)\ra R^{\cdot}(V_1\otimes V_2)$$
inducing the identity map on $V_1\otimes V_2$ and
satisfying the natural associativity condition
(e.g. one can take Cech complexes of acyclic covering or in the case $k=\C$
Dolbeault complexes). 
Then we can define a dg-category whose objects are vector bundles
with $\Hom(V_1,V_2)=R^{\cdot}(V_1^{\vee}\otimes V_2)$.
By homotopic invariance of the notion of $A_{\infty}$-algebra
there exists an equivalent
$A_{\infty}$-category structure on $\Vect(M)$
with morphisms $\Hom^*(V_1,V_2)=\oplus_i\Ext^i(V_1,V_2)$ which has $m_1=0$
(see \cite{GS},\cite{GLS},\cite{Kad},\cite{Markl},\cite{Merk}). 

We will use a particular representative of this class of $A_{\infty}$-structures
in the case when $M$ is a compact complex manifold equipped with a hermitian metric.
This $A_{\infty}$-structure appears naturally on the category $\Vect^h(M)$
of holomorphic vector bundles equipped with a hermitian metric.
Starting with the dg-category given by Dolbeault complexes one can use
metrics to write an explicit formula for higher compositions involving some
Hodge theory operators (see \cite{P-ainf}). We will denote this $A_{\infty}$-structure
by $m^H=(m^H_n)$.
Note that since $m_2^H$ is the standard composition (while $m_1=0$)
the choices of hermitian metrics on bundles are not really important. 
More precisely, by the standard argument in the homotopy theory
the objects $(V,h)$ and $(V,h')$, where $h$ and $h'$ are different metrics on the
same bundle $V$, are equivalent objects of this $A_{\infty}$-category. 
By the definition (that appeared in \cite{K-t}), this means that there exists
an $A_{\infty}$-functor from the category with two isomorphic objects
$O_1\simeq O_2$
and no other non-trivial morphisms to our $A_{\infty}$-category, that sends
$O_1$ to $(V,h)$ and $O_2$ to $(V,h')$. 

Assume in addition that $\om_M$ is trivialized. Then the Serre duality gives a
non-degenerate pairing
$$\Hom^*(V_1,V_2)\otimes\Hom^*(V_2,V_1)\ra\C.$$
The main feature of our particular choice of an $A_{\infty}$-structure
is the cyclic symmetry (\ref{cyclic}) of $m^H_n$ with respect to the Serre duality
(see \cite{P-ainf}).  
Note also that the higher products $m^H_n$
are compatible with Massey products when the latter are well-defined.

To define the second $A_{\infty}$-structure on $\Vect{E}$ (or rather, transversal
$A_{\infty}$-structure) 
let us recall the definition of the Fukaya $A_{\infty}$-category of the torus 
$T=\R^2/\Z^2$ with the (complexified) symplectic form 
$-2\pi i\tau dx\wedge dy$ where $\tau$ is an element of the upper half-plane.
We give here a very concrete version of the general definition which can be found
in \cite{F1},\cite{K},\cite{PZ}.
The objects of this category are pairs $(\ov{L},A)$
where  $\ov{L}=p(L)$ is the image of a non-vertical 
line $L$ with rational slope under the natural
projection $p:\R^2\ra T$ (a {\it geodesic circle}), 
$A:V\ra V$ is an operator on a finite dimensional complex vector space $V$
with real eigenvalues
\footnote{The slight difference with \cite{PZ} is that
we don't attach to $L$ an integer and don't consider vertical lines.}.
We call a pair of objects $(\ov{L_1},A_1)$ and $(\ov{L_2},A_2)$ transversal if
$\ov{L_1}$ and $\ov{L_2}$ are different. 
For such a pair the morphism space is 
$$\Hom((\ov{L_1},A_1),(\ov{L_2},A_2))=\Hom(V_1,V_2)\otimes\Hom(\ov{L_1},\ov{L_2})$$
where
$$\Hom(\ov{L_1},\ov{L_2})=\oplus_{P\in\ov{L_1}\cap\ov{L_2}}\C[P]$$ 
($[P]$ is a basis vector attached to a point $P$).
Note that there is a natural pairing
\begin{equation}\label{pairF}
\Hom((\ov{L_1},A_1),(\ov{L_2},A_2))\otimes\Hom((\ov{L_2},A_2),(\ov{L_1},A_1))
\ra\C
\end{equation}
induced by the natural duality between $\Hom(V_1,V_2)$ and $\Hom(V_2,V_1)$
and by the self-duality of $\Hom(\ov{L_1},\ov{L_2})=\Hom(\ov{L_2},\ov{L_1})$
(such that the basis $([P])$ is autodual).
Let $\la_i$ be the slope of the line $L_i$ ($i=1,2$). Then 
$\Hom((\ov{L_1},A_1),(\ov{L_2},A_2))\neq 0$ only if $\la_1\neq\la_2$.  
This space has grading $0$ if $\la_1<\la_2$ and grading $1$ if 
$\la_1>\la_2$.
By definition the differential $m_1$ is zero. The compositions $m_k$ 
for $k\ge 2$ are defined as follows.
Let $(\ov{L_i},A_i)$, $i=0,\ldots,k$
be objects of the Fukaya categories such that the corresponding
circles are pairwise different. Below it will be convenient
to identify the set of indices $[0,k]$ with $\Z/(k+1)\Z$.
For every $i\in\Z/(k+1)\Z$ let
$d_i\in [0,1]$ be the grading of $\Hom((\ov{L_i},A_i),(L_{i+1},A_{i+1}))$. The
composition
$$m_k^F:\Hom((\ov{L_0},A_0),(\ov{L_1},A_1))\otimes\ldots\otimes
\Hom((\ov{L_{k-1}},A_{k-1}),(\ov{L_k},A_k))\ra\Hom((\ov{L_0},A_0),(\ov{L_k},A_k))$$
is non-zero only if $\sum_{i=0}^{k}d_i=k-1$. 
Let $P_{i,i+1}$ be some intersection points of $\ov{L_i}$
and $\ov{L_{i+1}}$ ($i=0,\ldots,k-1$). 
For every $i=0,\ldots,k-1$ let $M_{i,i+1}$ be an element of 
$\Hom(V_i,V_{i+1})$. Then
\begin{align*}
&m_k^F(M_{0,1}[P_{0,1}],M_{1,2}[P_{1,2}],\ldots,M_{k-1,k}[P_{k-1,k}])=\\
&\sum_{P_{0,k},\Delta}\pm
\exp(2\pi i\tau\cdot\int_{\Delta}dx\wedge dy)
\exp(2\pi i(x(p_k)-x(p_{k-1}))A_k)\circ M_{k-1,k}\circ
\exp(2\pi i(x(p_{k-1})-x(p_{k-2}))A_{k-1})\\
&\ldots\circ M_{1,2}\circ\exp(2\pi i(x(p_1)-x(p_0))A_1)\circ M_{0,1}\circ
\exp(2\pi i(x(p_0)-x(p_k))A_0)
[P_{0,k}]
\end{align*}
where the sum is taken over points of intersections $P_{0,k}$ of
$\ov{L_0}$ with $\ov{L_k}$ and over all $(k+1)$-gons $\Delta$ 
(considered up to translation by $\Z^2$) with vertices
$p_i\equiv P_{i,i+1} \mod\Z^2$, $i\in\Z/(k+1)\Z$,
such that the edge $[p_{i-1},p_i]$ belongs to $p^{-1}(\ov{L_i})$
(the restriction on degrees $d_i$ implies
that $\Delta$ is convex). We also require
that the path formed by the edges $[p_0,p_1], [p_1,p_2], \ldots, [p_k,p_0]$
goes in the clockwise direction.
The sign in the RHS is given by the following rule.
If $k$ is even then all signs are ``plus''.
If $k$ odd then the sign is equal to the sign of $x(p_0)-x(p_k)$
(recall that we do not allow vertical lines).

It is not difficult to check that $m^F=(m^F_k)$ is a (transversal) 
cyclic $A_{\infty}$-category with respect to the pairing (\ref{pairF}). 
This $A_{\infty}$-structure
is strictly compatible with decomposition of an operator $A$ into
a direct sum of operators.
The main theorem of \cite{PZ} identifies the corresponding usual
category given by $m_2^F$ with a full subcategory of $\Vect(E)$ where
$E=\C/\Z+\Z\tau$ (which contains all indecomposable bundles). 
In order to get all vector bundles one has to modify
the Fukaya category by adding formally direct sums. We extend
the $A_{\infty}$-structure to this larger category using
strict compatibility with direct sums.
Thus, we get a transversal $A_{\infty}$-structure (which we still denote $m^F$) on 
$\Vect(E)$. The construction of \cite{PZ} identifies the pairing (\ref{pairF})
with the Serre duality (for some trivialization of $\om_E$), so the obtained 
$A_{\infty}$-structure is cyclic with respect to it.

We will remind some details of the correspondence between vector bundles on $E$
and objects of the Fukaya category later. Let us only mention here that
the slope of a line corresponding
to an indecomposable bundle $V$ is equal to the slope of $V$ (the ratio
of the degree and the rank). Stable bundles correspond to objects $(\ov{L},A)$
where $A\in\R$ is a real number 
(considered as an operator on a one-dimensional space).

\section{Transversal $A_{\infty}$-structures on the category
of line bundles over an elliptic curve}

\subsection{Transversality and admissibility}

Let $E$ be an elliptic curve over a field $k$.
Let $\LL$ be the full subcategory in $\Vect(E)$ consisting of line bundles. 
One can consider extensions of the (strictly
associative) composition $m_2$ on $\LL$ to $A_{\infty}$-structures.
The following definition gives some natural restrictions one can impose
on such an extension.
Let us fix a trivialization of $\om_E$. Then the Serre duality gives
a non-degenerate pairing
$$\Hom^*(V_1,V_2)\otimes\Hom^*(V_2,V_1)\ra k.$$

\begin{defi} Let us call a cyclic (with respect
to the Serre duality) $A_{\infty}$-structure $m$ on the category
$\LL$ {\it admissible} if $m_1=0$, $m_2$ is the standard
composition, and the functor of tensor multiplication by a line bundle is strictly
compatible with $m$. 
\end{defi}

Note that if $m_1=0$ then for any $A_{\infty}$-structure $m'$ which is
homotopic to $m$ one has $m'_2=m_2$. So it makes sense to try to
classify admissible $A_{\infty}$-structures on $\LL$ up
to cyclic homotopy, strictly compatible with
tensor multiplication by any line bundle. We refer to such homotopies as
{\it admissible} ones.

We also define an admissible {\it transversal} $A_{\infty}$-structure
on $\LL$ by similar restrictions provided that we have some notion of 
transversality for pairs of line bundles. We assume that such a notion
is given and that it has the following properties:

\noindent 
(i) $(L,M)$ is transversal if and only if $(M,L)$ is transversal;

\noindent
(ii) $(L,M)$ is transversal if and only if $(L^{-1},M^{-1})$ is transversal;

\noindent
(iii) $(L_1,L_2)$ is transversal if and only if $(L_1M,L_2M)$ is transversal;

\noindent
(iv) for every finite collection of line bundle $(L_1,\ldots,L_n)$
and every integer $d$ there exists an infinite number of
isomorphism classes of line bundles $L$ of degree $d$
such that $L$ and $L^2$ are transversal to all $L_i$;

\noindent
(v) if $(L,M)$ is transversal then $L\not\simeq M$.

For example, assume that $E(k)$ is infinite (this is necessary for
the property (iv) to hold). Then
one can call $(L,M)$ transversal if $L\not\simeq M$.
Another example arises from the correspondence between line bundles and objects
of Fukaya category defined in \cite{PZ}. 
In this example the complex parameter
describing an isomorphism class of a line bundle splits into two real parameters:
one describes the position of the corresponding geodesic circle and another
specifies the connection on it. Then the pair $(L,M)$ is transversal if the first
real parameter takes different values at $L$ and $M$ (see section \ref{conFuk}
for details).

The data of an admissible transversal $A_{\infty}$-structure are encoded in
the sequence of maps
$$m_n:H^{i_1}(L_1)H^{i_2}(L_2)\ldots H^{i_n}(L_n)\ra 
H^{i_1+i_2+\ldots+i_n+2-n}(L_1L_2\ldots L_n)$$
for line bundles $(L_i)$ such that the collection $(\O,L_1,L_1L_2,\ldots,
L_1\ldots L_n)$ is transversal. 

\subsection{Construction of the homotopy}

\begin{thm}\label{line-trans} 
Let $m$ and $m'$ be admissible transversal
$A_{\infty}$-structures
on the category of line bundles on $E$. Assume that
for every triple of line bundles
$(L_1,M,L_2)$ where $\deg(L_1)=\deg(L_2)=1$, $\deg(M)=-1$, such that
$(\O,L_1,L_1M,L_1ML_2)$ is transversal, the maps 
\begin{equation}\label{m3basic-l}
H^0(L_1)\otimes H^1(M)\otimes H^0(L_2)\ra H^0(L_1L_2M)
\end{equation}
given by $m_3$ and $m'_3$ coincide.
Then there exist a unique admissible homotopy between $m$ and $m'$.
\end{thm}

The following lemma is the main ingredient of the proof.

\begin{lem}\label{mainlem} 
Let $L$ be a line bundle of degree $\ge 3$ on $E$,
$S\subset\Pic(E)$ be a subset
such that for every $d$ and every isomorphism classes
$[L_1],[L_2]\in\Pic(E)$
there exists an infinite number of $[M]\in S$ such that
$\deg(M)=d$, $2[M]\in S$, $[L_1]+[M]\in S$ and $[L_2]-[M]\in S$.
Then the following sequence is exact:
\begin{equation}
\oplus_{L_1L_2L_3=L, [L_3]\in S, [L_2L_3]\in S}
H^0(L_1)H^0(L_2)H^0(L_3)\stackrel{\a}{\ra}
\oplus_{L_1L_2=L, [L_2]\in S} H^0(L_1) H^0(L_2) \stackrel{\b}{\ra} H^0(L)
\ra 0
\end{equation}
where $L_i$ denote line bundles of positive degrees, the map
$\a$ sends $s_1\otimes s_2\otimes s_3$ to
$s_1s_2\otimes s_3-s_1\otimes s_2s_3$, $\b$ sends $s_1\otimes s_2$ to
$s_1s_2$.
\end{lem}

\Pf . Clearly, $\b$ is surjective. Thus, it suffices to prove the following
statement. Assume that for every pair of line bundles of positive degree
$(L_1, L_2)$ such that $[L_2]\in S$ and
$L_1L_2\simeq L$ we are given a linear map
$$b_{L_1,L_2}:H^0(L_1)H^0(L_2)\ra k$$
such that for every triple $(L_1,L_2,L_3)$ such that $\deg(L_i)>0$,
$i=1,2,3$, $L_1L_2L_3\simeq L$, $[L_3]\in S$, $[L_2L_3]\in S$,
one has
$$b_{L_1,L_2L_3}(s_1,s_2s_3)=b_{L_1L_2,L_3}(s_1s_2,s_3)$$
where $s_i\in H^0(L_i)$, $i=1,2,3$. Then
there exists a functional $\phi$ on $H^0(L)$ such that
for every $(L_1,L_2)$ (with $[L_2]\in S$)
\begin{equation}\label{bL1L2}
b_{L_1,L_2}(s_1,s_2)=\phi(s_1s_2).
\end{equation}

We will consider separately several cases.

\noindent (i) $\deg(L)=3$.
Let $p_1,p_2\in E$ be a pair of distinct points such that
$[\O(p_i)]\in S$, $i=1,2$, and $[\O(p_1+p_2)]\in S$.
Let $s_{p_i}\in H^0(\O(p_i))$, $i=1,2$, be
non-zero sections.
Then we have the following exact sequence:
$$0\ra H^0(L(-p_1-p_2))\stackrel{\a'}{\ra} H^0(L(-p_1))\oplus H^0(L(-p_2))
\stackrel{\b'}{\ra} H^0(L)\ra 0$$
where $\a'(s)=(ss_{p_2},-ss_{p_1})$, $\b'(t_1,t_2)=t_1s_{p_1}+t_2s_{p_2}$.
Let us define a functional $\wt{\phi}$ on $H^0(L(-p_1))\oplus H^0(L(-p_2))$
by the formula
$$\wt{\phi}(t_1,t_2)=b_{L(-p_1),\O(p_1)}(t_1,s_{p_1})+
b_{L(-p_2),\O(p_2)}(t_2,s_{p_2}).$$
Note that $\wt{\phi}$ vanishes on the image of $\a'$.
Indeed, we have
$$b(ss_{p_2},s_{p_1})=b(s,s_{p_2}s_{p_1})=b(s,s_{p_1}s_{p_2})=
b(ss_{p_1},s_{p_2}).$$
Therefore, there exists a functional $\phi$ on $H^0(L)$ such that
$\wt{\phi}=\phi\circ\b'$. We are going to show that this functional is
the one we are looking for.

Let $L_1$ and $L_2$ be line bundles
of degrees $1$ and $2$ respectively such that $L_1L_2\simeq L$,
$[L_2]\in S$. Assume in addition that $L_2\not\simeq \O(p_1+p_2)$.
Then we claim that
$$b_{L_1,L_2}(s,t)=\phi(st)$$
for any $s\in H^0(L_1)$, $t\in H^0(L_2)$.
Indeed, the space $H^0(L_2)$ is a
direct sum of subspaces $H^0(L_2(-p_1))s_{p_1}$
and $H^0(L_2(-p_2))s_{p_2}$. Thus, it suffices to prove that
$b_{L_1,L_2}(s,t)=\phi(st)$ for $t$ in any of these subspaces.
For example, let $t=t's_{p_1}$, where $t'\in H^0(L_2(-p_1))$.
Then we have,
$$b(s,t)=b(s,t's_{p_1})=b(st',s_{p_1})=\phi(st's_{p_1})$$
as required.

Now we claim that if $M_1$ and $M_2$ are arbitrary line bundles
of degrees $2$ and $1$ respectively such that $M_1M_2\simeq L$
and $[M_2]\in S$, then $b_{M_1,M_2}(s,t)=\phi(st)$
for $s\in H^0(M_1)$, $t\in H^0(M_2)$. Indeed,
let us choose points $q_1,q_2\in E$ such that $M_1\not\simeq\O(q_1+q_2)$,
$M_2(q_i)\not\simeq\O(p_1+p_2)$ and $[M_2(q_i)]\in S$ for $i=1,2$.
Then $H^0(M_1)$ is a direct sum of $H^0(M_1(-q_1))H^0(\O(q_1))$ and
$H^0(M_1(-q_2))H^0(\O(q_2))$,
so we can assume that $s$ is in one of these subspaces. For example,
assume that $s=s's_{q_1}$ where $s_{q_1}\in H^0(\O(q_1))$,
$s'\in H^0(M_1(-q_1))$. Then we have
$$b(s,t)=b(s's_{q_1},t)=b(s',s_{q_1}t).$$
Applying the previous part of the proof to $L_1=M_1(-q_1)$,
$L_2=M_2(q_1)$ we obtain that
$$b(s',s_{q_1}t)=\phi(s's_{q_1}t)=\phi(st)$$
as required.

Finally, a similar argument shows that for arbitrary line bundles
$L_1$ and $L_2$ of degrees $1$ and $2$ one has
$b_{L_1,L_2}(s,t)=\phi(st)$ where $s\in H^0(L_1)$, $t\in H^0(L_2)$.

\noindent (ii) $\deg(L)=4$. Let us choose a pair of line bundles
$L_1$ and $L_2$ both of degree $2$ such that $L_1\not\simeq L_2$,
$[L_2]\in S$, and $L_1L_2\simeq L$. Then the product map
$$H^0(L_1)\otimes H^0(L_2)\ra H^0(L)$$
is an isomorphism, so we can define $\phi$ by setting
$$\phi(s_1s_2)=b_{L_1,L_2}(s_1,s_2)$$
where $s_1\in H^0(L_1)$, $s_2\in H^0(L_2)$. 

Let $L'_1$, $L'_2$ be line bundles of degree $2$ such that
$L'_1L'_2\simeq L$ and $[L'_2]\in S$. 
Let $p,q\in E$ be a pair of points such that $\O(p+q)\simeq L_1$,
$[L'_2(-q)]\in S$. 
Then we claim that the equality
\begin{equation}\label{bL'}
b_{L'_1,L'_2}(s,t)=\phi(st)
\end{equation}
holds whenever $s\in H^0(L'_1(-p))$, $t\in H^0(L'_2(-q))$.
Indeed, assume $s=s's_p$, $t=t's_q$ where
$s_p\in H^0(\O(p))$, $s_q\in H^0(\O(q))$, $s'\in H^0(L'_1(-p))$,
$t'\in H^0(L'_2(-q))$. Then
$$b(s,t)=
b(s's_p,t's_q)=b(s's_ps_q,t')=b(s_ps_q,s't')=\phi(s_ps_qs't')=\phi(st).$$

Now let $p_1,p_2\in E$ be a pair of distinct points 
such that $L'_1\simeq\O(p_1+p_2)$ and for $q_1,q_2\in E$
defined by $\O(p_1+q_1)\simeq\O(p_2+q_2)\simeq L_1$ one has
$[L'_2(-q_1)]\in S$, $[L'_2(-q_2)]\in S$.
Note that we have $L'_1(q_1+q_2)\simeq L_1^2\not\simeq L$
since $L_2\not\simeq L_1$. Hence, $\O(q_1+q_2)\not\simeq L'_2$
and $H^0(L'_2)$ has a basis $(t_1, t_2)$ such that 
$t_1$ vanishes at $q_1$ and $t_2$ vanishes at $q_2$.
Therefore, if $s$ is a section of $L'_1$ vanishing at $p_1$ and $p_2$
then by the previous part of the proof we have
$$b(s,t_i)=\phi(st_i)$$
for $i=1,2$. Repeating this argument for another pair of points 
$(p_1,p_2)$ as above we get a similar statement for a section of $L'_1$
linearly independent from $s$. Thus, we conclude that
(\ref{bL'}) holds for all $s\in H^0(L'_1)$, $t\in H^0(L'_2)$.

Now let $M_1$ be a line bundle of degree $1$, $M_2$ be a line bundle
of degree $3$ such that $M_1M_2\simeq L$, $[M_2]\in S$. 
Let $s_1\in H^0(M_1)$,
$s_2\in H^0(M_2)$, $p$ be a point in the divisor of $s_2$.
Then assuming that $[M_2(-p)]\in S$
we can write $s_2=s_ps'_2$ where $s_p\in H^0(\O(p))$,
$s'_2\in H^0(M_2(-p))$ and
$$b(s_1,s_2)=b(s_1,s_ps'_2)=b(s_1s_p,s'_2)=\phi(s_1s_ps'_2)=\phi(s_1s_2).$$
Since $H^0(M_2)$ is spanned by $H^0(M_2(-p))$ and $H^0(M_2(-p'))$
for two distinct points $p,p'\in E$, this proves that
$b_{M_1,M_2}(s_1,s_2)=\phi(s_1s_2)$ for all $s_1$ and $s_2$.
The case when $\deg(M_1)=3$, $\deg(M_2)=1$ is completely analogous.

\noindent (iii) $\deg(L)=d\ge 5$.
Let us fix a line bundle $L_2$ of degree $2$ such that $[L_2]\in S$
and $[L_2^2]\in S$.
Then there is an exact sequence
$$0\ra L_2^{-1}\ra H^0(L_2)\otimes\O\ra L_2\ra 0,$$
which induces for every line bundle $M$ of degree $\ge 3$ 
an exact sequence
$$0\ra H^0(ML_2^{-1})\ra H^0(M)H^0(L_2)\ra H^0(ML_2)\ra 0.$$
Let $L_1=LL_1^{-1}$. Consider the following diagram 
with exact rows
\begin{equation}
\begin{array}{ccccccc}
0\ra&\oplus H^0(\O(p))H^0(L_1L_2^{-1}(-p))&\ra
&\oplus H^0(\O(p))H^0(L_1(-p))H^0(L_2)&\stackrel{\a_2}\ra
&\oplus H^0(\O(p))H^0(L(-p))&\ra 0\\
&\ldar{\gamma}& &\ldar{\a_1} & &\ldar{}& \\
0\ra& H^0(L_1L_2^{-1})&\ra&H^0(L_1)H^0(L_2)&\stackrel{\b}\ra&H^0(L)&\ra0
\end{array}
\end{equation}
where the direct sums in first row are taken over all $p\in E$ such that
$[L(-p)]\in S$. Notice that $\gamma$ is surjective. Indeed, if $d\ge 6$
this is clear, while for $d=5$ we have to check that for the unique
point $p$ such that $\O(p)\simeq L_1L_2^{-1}$ one has $[L(-p)]\in S$.
But this follows from our assumptions of $L_2$, since 
$L(-p)\simeq L_2^2$ for such $p$.
He have $b_{L_1,L_2}\circ\a_1=\sum_p b_{\O(p),L(-p)}\circ\a_2$.
From this by an easy diagram chasing (using the surjectivity of $\gamma$)
we obtain that $b_{L_1,L_2}$
vanishes on the kernel of $\b$, hence, there exists a functional
$\phi$ on $H^0(L)$ such that $b_{L_1,L_2}=\phi\circ\b$.
It follows that for any $p\in E$ such that $[L(-p)]\in S$
one has 
$$b_{\O(p),L(-p)}(s,t)=\phi(st)$$ 
for $s\in H^0(\O(p))$, $t\in H^0(L(-p))$. Indeed,
we can assume that $t=t_1t_2$ with $t_1\in H^0(L_1(-p))$, $t_2\in L_2$,
in which case
$$b(s,t)=b(s,t_1t_2)=b(st_1,t_2)=\phi(st).$$
Now we can deduce (\ref{bL1L2}) in the general case using the
same argument as in
the end of case (ii).
\ed

\begin{rem} It is easy to see from the proof that our assumptions on
the set $S\in\Pic(E)$ can be weakened. Let us denote by $S_d$ the
subset of elements of $S$ of degree $d$. Then it suffices to require that:
(1) for any $[L]\in\Pic(E)$ and any $d$ one has $|S_d\cap (S+[L])|>4$,
$|S_d\cap ([L]-S)|>5$; (2) there exists $[L]\in S_2$ such that $2[L]\in S$.
\end{rem}

\noindent
{\it Proof of theorem \ref{line-trans}}. 
Let us prove the existence first.
Clearly we can replace $m$ by $m+\delta f$ where $f=(f_n,n\ge 2)$
is an admissible homotopy.
Therefore, we can argue by induction: for every $n\ge 3$,
assuming that $m_k=m'_k$ for $k<n$ we will construct an admissible homotopy
$f^n$ such that $f^n_k=0$ for $k<n-1$ and $(m+\delta f^n)_n=m'_n$
(this implies that $(m+\delta f^n)_k=m'_k$ for all $k\le n$). 

Using the cyclic symmetry we can reduce various types of
non-zero transversal $n$-tuple products to the following two types:

\noindent
(i)
$$m_n:H^0(L_1)H^1(M_1)\ldots H^1(M_i)H^0(L_2)H^1(M_{i+1})\ldots H^1(M_{n-2})\ra
H^0(L_1L_2M_1\ldots M_{n-2})$$
where $1\le i\le n-2$, 

\noindent
(ii)
$$m_n:H^0(L_1)\otimes H^0(L_2)\otimes H^1(M_1)\otimes\ldots\otimes
H^1(M_{n-2})\ra H^0(L_1L_2M_1\ldots M_{n-2}).$$

Let us call $w=\deg(L_1)+\deg(L_2)$ the {\it weight}
of the corresponding $n$-tuple product type. Note that
we have $w\ge 2$.
The first observation is that any $n$-tuple product of type (i)
of weight $>2$ can be expressed via
$k$-tuple products with $k<n$ and via $n$-tuple products of smaller
weight. Indeed, if $\deg(L_1)+\deg(L_2)>2$ then either $\deg(L_1)>1$
or $\deg(L_2)>1$. Assume for example that $\deg(L_1)>1$. Then
$H^0(L_1)$ is spanned by various products $s_ps$ where
$s_p\in \O(p)$, $s\in L_1(-p)$, a point $p\in E$ is such that the collection
$$(\O,\O(p), L_1,L_1M_1,\ldots,L_1M_1\ldots M_i,
L_1L_2M_1\ldots M_i, L_1L_2M_1\ldots M_{i+1},\ldots,
L_1L_2M_1\ldots M_{n-2})$$
is transversal.
Now for any collection of elements $e_j\in H^1(M_j)$, $j=1,\ldots,n-2$, 
$t\in H^0(L_2)$ and any $1\le j<n-2$ we have
\begin{align*}
&m_n(s_ps,e_1,\ldots,e_i,t,e_{i+1},\ldots,e_{n-2})=
m_n(s_p,se_1,\ldots,e_i,t,e_{i+1},\ldots,e_{n-2})\pm\\
&m_n(s_p,s,e_1,\ldots,e_it,e_{i+1},\ldots,e_{n-2})\pm
m_n(s_p,s,e_1,\ldots,e_i,te_{i+1},\ldots,e_{n-2})\pm\\ 
&s_p m_n(s,e_1,\ldots,e_i,t,e_{i+1},\ldots)+\ldots
\end{align*} 
where the unwritten terms contain only $m_k$ with $k<n$, while
the weights of three $n$-tuple products in the RHS are smaller than $w$.
If $j=n-2$ then there is an additional term
$m_n(s_p,s,e_1,\ldots,e_{n-2})t$ which doesn't affect our argument.
Similarly one considers the case when $\deg(L_2)>1$.
On the other hand, the only non-zero transversal products of type (i)
and weight $2$ are those of type (\ref{m3basic-l}).
As we'll see below this will allow us to restrict
our attention to products of type (ii).

To construct the homotopy $f^n$ we again apply induction. Namely,
assuming that $m_n=m'_n$ for all products (of types (i) and (ii))
of weight $<w$ (and $m_k=m'_k$ for $k<n$) we
will construct a homotopy $f^{n,w}$ such that
$(m+\delta f^{n,w})_n=m'_n$ for all products of weight $w$ and
such that the only non-zero component $f^{n,w}$ (other than $f^{n,w}_1=\id$)
reduces by cyclic symmetry to the following type:
$$f^{n,w}_{n-1}:H^0(L)\otimes H^1(M_1)\otimes\ldots\otimes H^1(M_{n-2})
\ra H^0(LM_1\ldots M_{n-2})$$
where $\deg(L)=w$. Note that $f^{n,w}$ is automatically cyclic. Indeed,
any non-zero value of $f^{n,w}$ is an element of $H^i(M)$ where the
degree of $M$ is either $-w$ or $d$,
such that $0<d<w$. On the other hand, by definition of Serre duality
$b(H^i(M),H^j(M'))=0$ unless $\deg(M)+\deg(M')=0$. It follows that one has
$$b(f^{n,w}_{n-1}(a_1,\ldots,a_{n-1}),f^{n,w}_{n-1}(a_n,\ldots,a_{2n-2}))=0,$$
so $f^{n,w}$ is cyclic. 
By the above observation it will be sufficient
to check the relation 
$(m+\delta f^{n,w})_n=m'_n$ only for products of type (ii) (and weight $w$).

Assume first that $w=2$. Then we necessarily have $n=3$.
Let us fix line bundles $L$ and $M$, $\deg(L)=2$, $\deg(M)=1$, such that
the triple $(\O,L,LM)$ is transversal.
We want to construct a map
$$f^{3,2}_2:H^0(L)\otimes H^1(M)\ra H^0(LM)$$
such that for every pair of line bundles
$L_1$, $L_2$ of degree $1$, where $L_1L_2\simeq L$
and the quadruple $(\O,L_1,L,LM)$ is transversal, the map
$$m'_3-m_3:H^0(L_1)H^0(L_2)H^1(M)\ra H^0(LM)$$
is a composition of the product map $H^0(L_1)H^0(L_2)\ra H^0(L)$
with $-f^{3,2}_2$.

Let us fix line bundles $M'$ and $L'$ such
that $\deg(M')=-2$, $\deg(L')=1$, $M'L'\simeq M$ and the quadruple
$(\O,L,LM',LM)$ is transversal.
Let $e\in H^1(M)$ be a non-zero element. Then $e=e's'$
for some $e'\in H^1(M')$, $s'\in H^0(L')$.
Now for every line bundles $L_1$ and $L_2$ such that
$L_1L_2\simeq L$, where the quintuple $(\O,L_1,L,LM',LM)$ is transversal, 
and every $s_1\in H^0(L_1)$ and $s_2\in H^0(L_2)$
we have
$$m_3(s_1,s_2,e)=m_3(s_1,s_2,e's')=m_3(s_1,s_2e',s')-m_3(s_1s_2,e',s')$$
and the similar equality holds for $m'_3$. Note that we have
$$m_3(s_1,s_2e',s')=m'_3(s_1,s_2e',s')$$
by the assumption of the theorem. Therefore,
\begin{equation}\label{hom1}
(m'_3-m_3)(s_1,s_2,e)=-(m'_3-m_3)(s_1s_2,e',s').
\end{equation}
Let us define the linear map
$$f_{e',s'}:H^0(L)\otimes H^1(M)\ra H^0(LM)$$
by the formula
$f_{e',s'}(s,e)=(m'_3-m_3)(s,e',s')$.
We claim that $f_{e',s'}$ doesn't depend on a choice
of $(M',L')$ and $e', s'$ such that $e's'=e$.
Indeed, $H^0(L)$ is generated by sections of the form $s=s_1s_2$ where
$s_1$ and $s_2$ are as above and the equality
(\ref{hom1}) shows that for such sections $f_{e',s'}(s,e)$
doesn't depend on $(e',s')$. Thus, we can set
$f^{3,2}_2=f_{e',s'}$. Now the same equality shows that
for any line bundles $L_1$, $L_2$ such
that $\deg(L_i)=1$, $L_1L_2\simeq L$ and the quadruple
$(\O,L_1,L,LM)$ is transversal one has
$$(m'_3-m_3)(s_1,s_2,e)=-f^{3,2}_2(s_1s_2,e).$$

Now assume that $w\ge 3$. 
Let us fix line bundles $M_1,\ldots, M_{n-2}$ and elements
$e_i\in H^1(M_i)$ for $i=1,\ldots,n-2$. Let us also fix a line
bundle $L$ of degree $w$, such that the collection
$$(\O,L,LM_1,\ldots,LM_1\ldots M_{n-2})$$ 
is transversal. Then for every
pair of line bundles $L_1$ and $L_2$ of positive degree such that
$L_1L_2\simeq L$ and the collection
$(\O,L_2,L_2M_1,\ldots,L_2M_1\ldots M_{n-2})$
is transversal, consider the map
$$b_{L_1,L_2}:H^0(L_1)H^0(L_2)\ra H^0(L_1L_2M_1\ldots M_{n-2}):
(s_1,s_2)\mapsto (m'_n-m_n)(s_1,s_2,e_1,\ldots,e_{n-2}).$$ 
We claim that these maps satisfy the condition
$$b(s_1s_2,s_3)=b(s_1,s_2s_3)$$
for any sections $s_i\in L_i$, $i=1,2,3$,
where $L_1L_2L_3\simeq L$, $\deg(L_i)>0$, the collection 
$$(\O,L_2,L_2L_3,L_2L_3M_1,\ldots,L_2L_3M_1\ldots M_{n-2})$$ 
is transversal.
Indeed, the constraint $\Ax_n$ implies that 
$$m_n(s_1s_2,s_3,e_1,\ldots,e_{n-2})-m_n(s_1,s_2s_3,e_1,\ldots,e_{n-2})$$
is a linear combination of terms either involving only $m_k$ with
$k<n$ or involving products $m_n$ of weight $<w$. The same is true for
$m'$, so our claim follows from the induction assumptions on $m$ and $m'$.
Therefore, we can
apply Lemma \ref{mainlem} to the line bundle $L$ and the set of isomorphism
classes 
$$S=\{ [M]:\ (\O,M,MM_1,\ldots,MM_1\ldots M_{n-2})
\ \text{is transversal} \}.$$
We conclude that there exists a linear map
$$f_{e_1,\ldots,e_{n-2}}:H^0(L)\ra H^0(LM_1\ldots M_{n-2})$$
satisfying
$$m'_n(s_1,s_2,e_1,\ldots,e_{n-2})-m_n(s_1,s_2,e_1,\ldots,e_{n-2})=
(-1)^n f_{e_1,\ldots,e_{n-2}}(s_1s_2,e_1,\ldots,e_{n-2}).$$
One can see from this defining property that the map
$$f^{n,w}_{n-1}:H^0(L)H^1(M_1)\ldots H^1(M_{n-2})\ra H^0(LM_1\ldots M_{n-2}):
s\otimes \otimes e_1\ldots\otimes e_{n-2}\ra f_{e_1,\ldots,e_{n-2}}(s)$$
is linear and gives the required homotopy.

The proof of uniqueness is also achieved by induction.
It suffices to check that an admissible transversal homotopy
$f=(f_n)$ from $m$ to $m$
such that $f_k=0$ for $2\le k<n$ has also $f_n=0$. By cyclic
symmetry it suffices to consider the maps
$$f_n:H^0(L)H^1(M_1)\ldots H^1(M_{n-1})\ra H^0(LM_1\ldots M_{n-1})$$
where $(\O,L,LM_1,\ldots,LM_1\ldots M_{n-1})$ is transversal.
Now we use the induction in degree of $L$. If $\deg(L)=1$ then such
a map is automatically zero. If $\deg(L)>1$ then it suffices to consider
elements of $H^0(L)$ of the form $ss_p$ where $s_p\in H^0(\O(p))$,
$s\in H^0(L(-p))$ (where $\O(p)$ is transversal to all the relevant bundles).
Then we can use the identity for $f_n$ and the induction
assumption to prove the desired vanishing.
\ed

\subsection{An identity between triple products}

Assume that we are given a transversal admissible $A_{\infty}$-structure
on the category of line bundles on $E$.
Let $(L_1,M,L_2)$ be a triple of line bundles such that
$\deg(L_1)=\deg(L_2)=n>0$, $\deg(M)=-n$, and the collection
$(\O,L_1,L_1M,L_1ML_2)$ is transversal. Then the triple
products 
$$m_3:H^0(L_1)H^1(M)H^0(L_2)\ra H^0(L_1ML_2)$$
are invariant under any homotopy. However, in theorem \ref{line-trans}
only such triple products with $n=1$ appear. The reason is
that one can express all triple products as above in terms of those with
$n=1$. This is done by induction in $n$ using the identity below.

Assume that $L_i=L'_iL''_i$ for $i=1,2$, where $\deg(L'_i)=n'$,
$\deg(L''_i)=n''$ for some positive integers $n'$, $n''$ such that
$n=n'+n''$. Assume also that the collection
$(\O,L'_1,L_1,L_1M,L_1ML'_2,L_1ML_2)$ is transversal.

\begin{prop} One has the following identity
$$m_3(s'_1s''_1,e,s'_2s''_2)=m_3(s'_1,s''_1e,s'_2)s''_2+
s'_1 m_3(s''_1,es'_2,s''_2)$$
where $s'_i\in H^0(L'_i)$, $s''_i\in H^0(L''_i)$, $e\in H^1(M)$.
\end{prop}

\Pf . Applying the $A_{\infty}$-constraint $\Ax_3$ we get
$$m_3(s'_1s''_1,e,s_2)=m_3(s'_1,s''_1e,s'_2s''_2)+
s'_1 m_3(s''_1,e,s'_2s''_2).$$
Applying $\Ax_3$ again we obtain the following expressions for the terms
in the RHS:
$$m_3(s'_1,s''_1e,s'_2s''_2)=m_3(s'_1,s''_1e,s'_2)s''_2+
s'_1 m_3(s''_1e,s'_2,s''_2),$$ 
$$m_3(s''_1,e,s'_2s''_2)=m_3(s''_1,es'_2,s''_2)-m_3(s''_1e,s'_2,s''_2).$$
Substituting these expressions in the above equality we get the result.
\ed

\section{Application to homological mirror symmetry}

\subsection{Adding unipotent bundles}\label{uni-sec}

By a unipotent bundle we mean a vector bundle which has a filtration by subbundles
such that the associated graded bundle is trivial.
Let $\LL\UU=\LL\UU(E)$ be the full subcategory in $\Vect(E)$ consisting of bundles
of the form $LU$, where $L$ is a line bundle, $U$ is a unipotent bundle.
Note that a decomposition of $LU$ into a tensor product of a line bundle
and a unipotent bundle is unique up to an isomorphism.

Assume that we are given a notion of transversality for
pairs of line bundles. We can extend it to the category $\LL\UU$ by
calling a pair $(LU, L'U')$ transversal if and only if $(L,L')$ is transversal.
Then we define an admissible transversal $A_{\infty}$-structure on $\LL\UU$
as a transversal $A_{\infty}$-structure on $\LL\UU$ which is cyclic
with respect to Serre duality and is
strictly compatible with tensor multiplication by a line bundle, has $m_1=0$
and $m_2$ equal to the standard product.

One defines a notion of admissible homotopy between admissible $A_{\infty}$-structures
on $\LL\UU$ similarly to the case of the category $\LL$.

The proof of the following theorem is very similar to that of
theorem \ref{line-trans} so we omit it.

\begin{thm}\label{unip} Let $m$ and $m'$ be admissible transversal
$A_{\infty}$-structures on the category $\LL\UU$. Assume that
for every triple of line bundles
$(L_1,M,L_2)$ such that $\deg(L_1)=\deg(L_2)=1$, $\deg(M)=-1$ and such that
$(\O,L_1,L_1M,L_1ML_2)$ is transversal, and for every quadruple of unipotent
bundles $U_0$, $U_1$, $U_2$ and $U_3$ the maps 
\begin{equation}\label{m3basic}
\Hom(U_0,L_1U_1)\otimes \Ext^1(L_1U_1,L_1MU_2)\otimes 
\Hom(L_1MU_2,L_1ML_2U_3)\ra \Hom(U_0,L_1ML_2U_3)
\end{equation}
given by $m_3$ and $m'_3$ coincide.
Then there exist a unique admissible homotopy between $m$ and $m'$.
\end{thm}

\subsection{Connection with the Fukaya category}\label{conFuk}

Let $\tau\in\C$ be an element in the upper half-plane,
$E=\C/\Z+\Z\tau$, be the corresponding elliptic curve.
Then as shown in \cite{PZ} the (usual) category $\LL\UU$
is equivalent to the subcategory in the Fukaya category (with
compositions $m_2^F$) consisting of objects $(\ov{L},\la\cdot\Id+N)$ where
$\ov{L}$ has an integer slope, $\la\in\R$, $N$ is a nilpotent operator.

Let $L(0)$ be the line bundle on $E$ such that the theta-function
$$\th(z)=\th(z,\tau)=\sum_{n\in\Z}\exp(\pi i \tau n^2+2\pi i n z)$$
is the pull-back of a section of $L$. So $L(0)\simeq\O_E(z_0)$
where $z_0=\frac{\tau+1}{2}\mod(\Z+\Z\tau)$.
For every $u\in\C$ let us denote $L(u)=t^*_u L(0)$, where $t_u:E\ra E$ is the
translation by $u$. Then
every line bundle of degree $n$ is isomorphic to a line bundle
of the form $L(0)^{\otimes(n-1)}\otimes L(u)$.

For a nilpotent operator $N:V\ra V$
we denote by $\VV_N$ the unipotent bundle on $E$, such that the sections
of $\VV_N$ correspond to $V$-valued functions on $\C$ satisfying the
quasi-periodicity equations
$f(z+1)=f(z)$, $f(z+\tau)=\exp(2\pi i N)f(z)$.
Then every unipotent bundle is isomorphic to a bundle of the form $\VV_N$.

The correspondence between bundles in $\LL\UU$ and objects of the Fukaya
category constructed in \cite{PZ} associates to the bundle
$\VV=L(0)^{\otimes(n-1)}\otimes L(u)\otimes\VV_N$ the object
$O=(\ov{L}, -u_1\Id+N)$, where $u=u_1+\tau u_2$, $u_i\in\R$,
$\ov{L}=\{(u_2+x,(n-1)u_2+nx), x\in\R/\Z\}$. 

This correspondence extends to a functor from $\LL\UU$ to the Fukaya
category (with $m_2^F$ as a composition) as follows.
Let $\VV'=L(0)^{\otimes(n'-1)}\otimes L(u')\otimes\VV_{N'}$ be
another bundle in $\LL\UU$, where $n'\in\Z$, $u'=u'_1+\tau u'_2\in\C$, 
$N':V'\ra V'$
is a nilpotent operator. Let $O'=(\ov{L'},-u'_1\Id+N')$ be the corresponding
object in the Fukaya category. Note that $O$ and $O'$ are transversal
if and only if either $n'\neq n$, or $n'=n$ and $u'_2-u_2\not\in\Z$.
In the latter case $\Hom(\VV,\VV')=\Hom(O,O')=0$ so we can assume that $n\neq n'$.
Assume first that $n<n'$.
Then $\Hom(O,O')=\Hom(V,V')\otimes\Hom(\ov{L},\ov{L'})$ has degree zero.
We can enumerate the points of intersection $\ov{L}\cap\ov{L'}$ by 
residues $k\in\Z/(n'-n)\Z$. Namely, this intersection consists of the points 
$$P_k=(\frac{k+u'_2-u_2}{n'-n},\frac{nk+nu'_2-n'u_2}{n'-n})$$
where $k\in\Z/(n'-n)\Z$.
On the other hand, we have
$$\Hom(\VV,\VV')=H^0(E,L(0)^{n'-n-1}\otimes L(u'-u)\otimes\VV_{N'-N^*})$$
where we consider $N^*$ and $N'$ as operators on $V^*\otimes V'$ (acting
trivially on one component).
Note that if $M$ is a line bundle on $E$ of the form 
$L(0)^{\otimes(m-1)}\otimes L(u)$ where $m\neq 0$ and
$N:V\ra V$ is a nilpotent operator then there is a natural isomorphism
between Dolbeault complexes of bundles $M\otimes V$ and $M\otimes_{\O}\VV_N$.  
Indeed, using the trivialization of the pull-backs of $M$ and $\VV_N$
to $\C$ we can
define the map from the Dolbeault complex of $M\otimes V$ to that of
$M\otimes_{\O}\VV_N$ by sending $\eta(z)$ to $\eta(z-N/m)$ where
$$(f\otimes v)(z-\frac{N}{m})=\exp(-\pa_z \frac{N}{m})(f)\cdot v$$
In particular, we can identify $\Hom(\VV,\VV')$ with 
the space $\Hom(V,V')\otimes H^0(E,L(0)^{n'-n-1}\otimes L(u'-u))$.
The space of global sections of the line bundle
$L(0)^{n'-n-1}\otimes L(u'-u)$ has a natural basis of theta functions
$$\th_k(z)=\sum_{m\in (n'-n)\Z+k}\exp(\frac{1}{n'-n}(\pi i\tau m^2+
2\pi i m((n'-n)z+u'-u)))$$
where $k\in\Z/(n'-n)\Z$.
Now we can identify $\Hom(\VV,\VV')$ with $\Hom(O,O')$ by sending
$T\otimes [P_k]$ (where $T\in\Hom(V,V')$) to 
$$\exp(\frac{1}{n'-n}(-\pi i\tau (u'_2-u_2)^2\Id+
2\pi i(u'_2-u_2)(N'-N^*-(u'_1-u_1)\Id)))\cdot T\otimes \th_k.$$

To construct similar identification in the case $n>n'$ 
we use Serre duality and its natural analogue on the Fukaya
category to reduce to the case considered above. 
As shown in \cite{PZ} this identification is compatible with compositions
$m_2$. Using it we can consider $m^F$ as a transversal $A_{\infty}$-structure on 
$\LL\UU$. Furthermore, it is easy to see that $m^F$ is admissible. The main
point is that the functor of tensoring with a line bundle on $\LL\UU$ corresponds
to an automorphism of the Fukaya category given by some symplectic automorphism of 
the torus. As we will see in section \ref{massey} the assumptions of the theorem
\ref{unip} are satisfied for the transversal
$A_{\infty}$-structures $m^F$ and $m^H$ on $\LL\UU$. Hence, they are homotopic.

The equivalence of $\LL\UU$ with a subcategory of the Fukaya category (with
$m_2^F$ as a composition) is extended to all bundles in \cite{PZ} using
the construction of vector bundles on $E$ as push-forwards of
objects in $\LL\UU$ under isogenies. Below we consider the corresponding
extension of equivalence between $A_{\infty}$-structures.

\subsection{Equivalence}

Let $\tau\in\C$ be an element in the upper half-plane,
$E=\C/\Z+\Z\tau$, be the corresponding elliptic curve.
We are going to prove that the two transversal $A_{\infty}$-structures 
$m^F$ and $m^H$ on
$\Vect(E)$ considered in section \ref{two} are equivalent.
More precisely, the definition of $m^H$ requires us to work with bundles equipped 
with hermitian metrics.
Since the different choices of metrics on a vector bundle really give equivalent 
objects of the $A_{\infty}$-category $(\Vect^h(E),m^H)$ 
we can restrict to some preferred class of hermitian metrics 
(which we'll define below).

Note that both $A_{\infty}$-structures are cyclic and
are strictly compatible with tensor multiplication by a hermitian line bundle
and with decompositions of bundles into orthogonal direct sums. 

For every positive integer
$r$ we consider the elliptic curve $E_r=\C/\Z+\Z r\tau$. Then we have
a natural isogeny $\pi^r:E_r\ra E$ of degree $r$
and for every $r|s$ an isogeny $\pi^s_r:E_s\ra E_r$ such that
$\pi^s=\pi^r\circ\pi^s_r$. 
We can consider two transversal
$A_{\infty}$-structures $m^F$ and $m^H$ on any of these
elliptic curves. An important observation is that
both $m^F$ and $m^H$ are strictly compatible with the functors
of pull-back and push-forward with respect to isogenies $\pi^s_r$ and $\pi^r$
(these functors extend naturally to vector bundles with metric).
For the structure $m^H$ this is clear while for $m^F$ this follows from the 
construction of equivalence in \cite{PZ}.

The idea of the proof is to use the decomposition of every bundle on elliptic curve
into a direct sum 
$V=\oplus V_iU_i$ where $(V_i)$ are pairwise non-isomorphic stable bundles,
$(U_i)$ are unipotent bundles. Then we want to use the fact that every stable
bundle of rank $r$ on $E$ is the push-forward of a line bundle on $E_r$. 
Since our $A_{\infty}$-structures are strictly compatible with
isogenies we can derive the desired homotopy from theorem \ref{unip}.
More precisely, we need a slight modification of this theorem for the category
of bundles with metrics: the assumption should be that the triple products
(\ref{m3basic}) given by two $A_{\infty}$-structures coincide for all choices
of metrics on the bundles in question.
To be able to apply this theorem in our case we will compute explicitly 
products $m_3^F$ and $m_3^H$ of the type (\ref{m3basic}) in section \ref{massey}
and will see that they
are equal (at this point it will be important to use a particular trivialization
of $\om_E$ on which the construction of equivalence in \cite{PZ} depends).
Then the uniqueness of the homotopies constructed in theorem \ref{unip}
will imply that these homotopies are compatible with isogenies, hence, descend
to a homotopy on the category $\Vect(E)$.

Let us call a vector bundle on $E$ {\it almost stable} if it has form
$V\otimes U$ where $V$ is a stable bundle, $U$ is a unipotent bundle
(thus, every bundle on $E$ is a direct sum of almost stable bundles). 
Let $r$ be the rank of $V$. Then $V=\pi^r_*(L)$ for some line bundle
$L$ on $E_r$. Hence $V\otimes U=\pi^r_*(L\otimes(\pi^r)^*U)$.
We call a hermitian metric on $V$ {\it preferred} if
it comes from a metric on $L\otimes(\pi^r)^*U$.
Let $(V_1,\ldots,V_{n+1})$ be a collection of almost stable bundles on $E$ 
equipped with preferred metrics.
From the strict compatibility of our $A_{\infty}$-structures with isogenies 
we have the following commutative diagram
\begin{equation}\label{mn-compat}
\begin{array}{ccccccc}
\Hom^*(V_1,V_2)\ldots\Hom^*(V_n,V_{n+1})&\lrar{m_n}&\Hom^*(V_1,V_{n+1})\\
\ldar{(\pi^r)^*} & & \ldar{(\pi^r)^*}\\
\Hom^*((\pi^r)^*V_1,(\pi^r)^*V_2)\ldots\Hom((\pi^r)^*V_n,(\pi^r)^*V_{n+1})&\lrar{m_n}
&\Hom^*((\pi^r)^*V_1,(\pi^r)^*V_{n+1})
\end{array}
\end{equation}
where $m=m^F$ or $m=m^H$.
Now if $r$ is divisible by ranks of all bundles $V_i$ then 
$(\pi^r)^*(V_i)$ is an orthogonal direct
sum of bundles of the form $LU$ where $L$ is a line bundle, $U$ is a unipotent bundle.

We will check in section \ref{massey} that the conditions of theorem
\ref{unip} are satisfied for $m^F$ and $m^H$. Therefore, we get a unique admissible
homotopy $f^r$ between these structures on $\LL\UU(E^r)$ for every $r$.
We can extend this homotopy to orthogonal direct sums of bundles in $\LL\UU(E^r)$
in an obvious way. Note that for every isogeny of elliptic curves
$\pi:E'\ra E''$ we have a canonical splitting of the natural embedding
$\O_{E''}\ra\pi_*\O_{E'}$, hence for every pair of bundles $(V_1,V_2)$ on $E''$
we get a canonical splitting
$$T(\pi):\Hom(\pi^*V_1,\pi^*V_2)\ra\Hom(V_1,V_2)$$
of the natural embedding $\pi^*:\Hom(V_1,V_2)\ra\Hom(\pi^*V_1,\pi^*V_2)$.
Now we claim that the homotopies $f^r$ and $f^s$ where $r|s$
are compatible in the following way:
for any bundles $W_1=L_1U_1,\ldots,W_{n+1}=L_{n+1}U_{n+1}$ in $\LL\UU(E^r)$
one has the commutative diagram
\begin{equation}\label{f-comp}
\begin{array}{ccccccc}
\Hom^*(W_1,W_2)\ldots\Hom^*(W_n,W_{n+1})&\lrar{f^r_n}&
\Hom^*(W_1,W_{n+1})\\
\ldar{(\pi^s_r)^*} & & \luar{T(\pi^s_r)}\\
\Hom^*((\pi^s_r)^*(W_1),(\pi^s_r)^*(W_2))\ldots\Hom((\pi^s_r)^*(W_n),
(\pi^s_r)^*(W_{n+1}))
&\lrar{m_n}&\Hom^*((\pi^s_r)^*(W_1),(\pi^s_r)^*(W_{n+1}))
\end{array}
\end{equation}
Indeed, the compatibility of $m^F$ and $m^H$ with the isogeny $\pi^s_r$ implies
that $T(\pi^s_r)\circ f^s\circ(\pi^s_r)^*$ is an admissible homotopy between
$m^F$ and $m^H$ on $\LL\UU(E^r)$, hence it coincides with $f^r$.

Now we define the homotopy $f$ between $m^F$ and $m^H$ on the category
of almost stable vector bundles on $E$ with preferred metrics using commutativity
of diagrams of the type (\ref{mn-compat}). 
Namely, choosing $r$ which is divisible by all
ranks of bundles $V_i$ we define the map
$$f_n:\Hom(V_1,V_2)\ldots\Hom(V_n,V_{n+1})\ra\Hom(V_1,V_{n+1})$$ 
by the formula
$f_n=T(\pi^r)\circ f^r_n\circ (\pi^r)^*$.
The compatibility (\ref{f-comp}) ensures that this definition doesn't depend
on a choice $r$.
Now to check that $f$ is indeed a homotopy from $m^F$ to $m^H$ we choose
$r$ divisible by ranks of all the bundles involved and use the commutativity
of (\ref{mn-compat}).

Since every bundle $V$ on $E$ is a direct sum of almost stable bundles we have
a class of preferred metrics on $V$ coming from preferred metrics on almost
stable bundles (so that the direct sum becomes orthogonal). We can extend
the homotopy $f$ to all bundles with preferred metrics in a natural way. 

\subsection{Massey products}\label{massey}
It remains to compute explicitly the products $m_3^F$ and $m_3^H$ of the type 
(\ref{m3basic}).
Let us trivialize $\om_E$ in such a way that the Serre duality induces
the pairing 
$$b:\Hom(V_1,V_2)\otimes\Ext^1(V_2,V_1)\ra\C$$
given by the formula
$$b(f,gd\ov{z})=\int_{E} dz\we \Tr(f\circ g d\ov{z})$$
where $f\in\Hom(V_1,V_2)$, $gd\ov{z}\in\Om^{0,1}(\Hom(V_2,V_1))$.

First let us compute $m_3^H$. We start with the case when
all $U_i$ are trivial of rank $1$.
Then we have to compute the product
$$m_3^H:H^0(L_1)H^1(M)H^0(L_2)\ra H^0(L_1ML_2)$$
where $L_1M\not\simeq\O$, $L_2M\not\simeq\O$. 
Using a translation on $E$ we can assume without loss of generality that
$M=L(0)^{-1}$. Let $L_1=L(t)$, $L_2=L(u)$ where $t,u\in\C$.
Let $z_1$ and $z_2$ be the real components of the complex variable $z$
defined by the equality $z=z_1+\tau z_2$. The transversality condition
means that $t_2,u_2\not\in\Z$. We will compute the above
product under the weaker assumption $t,u\not\in\Z+\Z\tau$.
It is easy to check that the $(0,1)$-form with values in $L(0)^{-1}$
$$\a(z)=\frac{i}{\sqrt{2\Im(\tau)}}\ov{\theta(z)}
\exp(-2\pi \Im(\tau)(z_2^2))d\ov{z}$$
is a representative of the class in $H^1(L(0)^{-1})$ dual to the class
in $H^0(L(0))$ given by $\th(z)$.
Now for every $u\in\C$, such that $u\not\in\Z+\Z\tau$ there exists
a unique section $h(z,u)$ of $L(0)^{-1}L(u)$ such that
$$\th(z+u)\a(z)=\ov{\pa}h(z,u)$$
where $\ov{\pa}=\ov{\pa}_z$. Indeed, this follows from the fact
that all the cohomologies of $L(0)^{-1}L(u)$ vanish.
One can write an explicit formula
for $h(z,u)$(see \cite{P-ainf}):
$$h(z,u)=-\frac{1}{2\pi i}\sum_{m,n\in\Z}(-1)^{mn}\frac
{\exp(-\frac{\pi}{2\Im(\tau)}(|\ga|^2+2\ov{\ga}u+u^2)+2\pi i(mz_1+(n-u)z_2))}{\ga+u}$$
where $\ga=m\tau-n$.
Now we have
$$m_3^H(\th(z+t),\a,\th(z+u))=h(z,t)\th(z+u)-h(z,u)\th(z+t).$$
As a function of $z$ up to a constant factor this should be equal to
$\th(z+u+v)$, so we have
\begin{equation}\label{Htu}
h(z,t)\th(z+u)-h(z,u)\th(z+t)=H(t,u)\th(z+t+u)
\end{equation}
for some meromorphic function $H$.
We have $H(t,u)=-H(u,t)$. Also it is easy to see that the function $H(t,u)$
satisfies the following quasi-periodicity equations:
$$H(t+1,u)=H(t,u),$$
$$H(t+\tau,u)=\exp(2\pi i u)H(t,u).$$
The only poles of $H(t,u)$ are poles of order $1$ along the divisors
$t=\ga$ and $u=\ga$ where $\ga\in\Z+\Z\tau$. It follows that $H(t,u)$
is equal up to a constant to the function
$$F(t,u)=\frac{\th'(\frac{\tau+1}{2})\th(t-u+\frac{\tau+1}{2})}
{2\pi i\th(t+\frac{\tau+1}{2})\th(-u+\frac{\tau+1}{2})}.$$
Furthermore, comparing the residues at $t=0$ we conclude that
$H(t,u)=-F(t,u)$.

Now let us compute the product
$$m_3^H:H^0(L_1U_0^{\vee}U_1)H^1(MU_1^{\vee}U_2)H^0(L_2U_2^{\vee}U_3)\ra
H^0(L_1ML_2U_0^{\vee}U_3)$$
where $U_i$ are unipotent bundles. As before we can take
$M=L(0)^{-1}$, $L_1=L(t)$, $L_2=L(u)$.
Let $U_i=\VV_{N_i}$ where $N_i:V_i\ra V_i$
are nilpotent operators. Then $U_i^*U_{i+1}\simeq\VV_{N_{i+1}-N_i^*}$
where $N_{i+1}-N_i^*$ is an operator on $V_i^*V_{i+1}$.
As in section \ref{conFuk} we use the
isomorphisms between the Dolbeault complexes
of bundles $LV$ and $L\VV_N$, 
where $L$ is one of line bundles of degree $1$ above, $N:V\ra V$
is the corresponding nilpotent operator, sending
$\eta(z)$ to $\eta(z-N)$.
Similarly, we have an isomorphism between the Dolbeault complexes of
$L(0)^{-1}V_1^*V_2$ and $L(0)^{-1}\VV_{N_2-N_1^*}$
given by $\eta(z)\mapsto \eta(z+N_2-N_1^*)$.
Let $v_{i,i+1}\in V_i^*\otimes V_{i+1}$ be
some elements. Then we have
\begin{align*}
&(\a(z+N_2-N_1^*)v_{1,2})\circ(\th(z+t-N_1+N_0^*))v_{0,1})=\\
&\Tr_{V_1}(\ov{\pa}h(z+N_2-N_1^*,t-N_2+N_1^*-N_1+N_0^*))
v_{0,1}v_{1,2})=\\
&\Tr_{V_1}(\ov{\pa}h(z+N_2-N_1^*,t-N_2+N_0^*)v_{0,1}
v_{1,2}),
\end{align*}
since we can replace $N_1^*$ by $N_1$ under the sign of $\Tr_{V_1}$.
Similarly, we get
$$(\th(z+u-N_3+N_2^*))v_{2,3})\circ(\a(z+N_2-N_1^*)v_{1,2})=
\Tr_{V_2}(\ov{\pa}h(z+N_2-N_1^*,u-N_3+N_1^*)v_{1,2}v_{2,3}).$$
Hence,
\begin{align*}
&m_3^H(\th(z+t-N_1+N_0^*))v_{0,1},\a(z+N_2-N_1^*)v_{1,2},
\th(z+u-N_3+N_2^*)v_{2,3})=\\
&\Tr_{V_1V_2}(
(\th(z+u-N_3+N_2^*)h(z+N_2-N_1^*,t-N_2+N_0^*)-\\
&h(z+N_2-N_1^*,u-N_3+N_1^*)\th(z+t-N_1+N_0^*))v_{0,1}v_{1,2}v_{2,3}).
\end{align*}
Making a substitution $z\mapsto z+N_2-N_1^*$,
$t\mapsto t-N_2+N_0^*$, $u\mapsto u-N_3+N_1^*$ in the identity
(\ref{Htu}) and using the equality $H=-F$
we can rewrite the above formula as follows:
\begin{equation}\label{m3H-main}
\begin{array}{c}
m_3^H(\th(z+t-N_1+N_0^*)v_{0,1},\a(z+N_2-N_1^*)v_{1,2},
\th(z+u-N_3+N_2^*)v_{2,3})=\\
\Tr_{V_1V_2}(F(t-N_2+N_0^*),u-N_3+N_1^*))
\th(z+t+u-N_3+N_0^*)v_{0,1}v_{1,2}v_{2,3}).
\end{array}
\end{equation}

Now let us compute the corresponding product $m_3^F$.
The objects of the Fukaya category corresponding to
our four bundles $U_0=\VV_{N_0}$, $L_1U_1=L(t)\VV_{N_1}$,
$L_1MU_2=L(0)^{-1}L(t)\VV_{N_2}$ and $L_1ML_2U_3=L(t+u)\VV_{N_3}$
are $((x,0),N_0)$, $((x+t_2,x),-t_1+N_1)$,
$((x,-t_2),-t_1+N_2)$ and $((x+t_2+u_2,x),-t_1-u_1+N_3)$,
where $t=t_1+\tau t_2$, $u=u_1+\tau u_2$, $t_2,u_2\not\in\Z$.
Note that any two of these circles either don't intersect or intersect
at a unique point. So we can identify morphisms between these objects
with spaces $\Hom(V_0,V_1)$, $\Hom(V_1,V_2)$, etc.
Now we have
\begin{align*}
&-m_3^F(v_{0,1}[P_{0,1}],v_{1,2}[P_{1,2}],v_{2,3}[P_{2,3}])=
\Tr_{V_1V_2}
\sum_{(m,n)\in\Z^2,(m-t_2)(n+u_2)>0}\sign(m-t_2)\\
&\exp(2\pi i\tau(m-t_2)(n+u_2)+
2\pi i(m-t_2)(-t_1+N_2-N_0^*)+2\pi i(n+u_2)(u_1-N_3+N_1^*))
v_{0,1}v_{1,2}v_{2,3})[P_{0,3}]\\
&=\Tr_{V_1V_2}(\sum\sign(m-t_2)\exp(2\pi i\tau mn+
2\pi i m(u-N_3+N_1^*)+2\pi i n(-t+N_2-N_0^*))\cdot C\cdot
v_{0,1}v_{1,2}v_{2,3})
\end{align*}
where $C=\exp(-2\pi i\tau t_2u_2-2\pi i t_2(u_1-N_3+N_1^*)+2\pi i
u_2(-t_1+N_2-N_0^*))$.
At this point we need the following identity
(which essentially coincides with the formula (2.3.4) of \cite{P-Mas}): 
$$\sum_{(m,n)\in\Z^2,(m-t_2)(n+u_2)>0}\sign(m-t_2)
\exp(2\pi i\tau mn+2\pi i (mu-nt))=F(t,u)$$
for arbitrary $t=t_1+\tau t_2$, $u=u_1+\tau u_2$ such that
$t_2,u_2\in\Z$.
This identity which is due to Kronecker can be proven as follows: 
first, one has to check that the
left hand side extends to a meromorphic function of $u$ and $t$ with
poles at the lattice points, then one has to compare its quasi-periodicity properties
and residues at poles with those of $F$. 
Hence, we get
\begin{equation}\label{m3F-main}
m_3^F(v_{0,1}[P_{0,1}],v_{1,2}[P_{1,2}],v_{2,3}[P_{2,3}])=
-\Tr_{V_1V_2} (F(t-N_2+N_0^*,u-N_3+N_1^*)\cdot C\cdot
v_{0,1}v_{1,2}v_{2,3})[P_{0,3}].
\end{equation}
Now it easy to see that the exponential factors involved in the
identification of morphisms in $\LL\UU$ with morphisms in the Fukaya
category (see section \ref{conFuk}) kill the factor $C$
and we get $m_3^H=m_3^F$ on the products of the type
(\ref{m3basic}).

\end{document}